\documentclass[a4paper]{jpconf}
\usepackage{graphicx}

\usepackage{amsmath,amssymb,cite}

\def\beq{\begin{equation}}
\def\eqn#1{\beq\label{#1}}
\def\eeq{\end{equation}}

\def\bb {\begin {eqnarray}}
\def\eqnn#1{\bb\label{#1}}
\def\ee {\end {eqnarray}}

\def\bbz{\mathbb{Z}}
\def\bbc{\mathbb{C}}

\def\bbr{\mathbb{R}}
\def\bbn{\mathbb{N}}

\def\nn{\nonumber}
\def\nt{\noindent}

\def\np{\vfill\eject}
 \def\nd{\end{document}}

\def\ftha{{\textstyle{15\over2}}}

\def\llr{\longrightarrow}
\def\({\left(}
\def\){\right)}
\def\eps{\epsilon}
 
\def\lra{\longrightarrow}

\def\lg{\langle} \def\rg{\rangle} 

\def\ha{{\textstyle{1\over2}}}

  \def\tV{{\tilde V}}

\def\r{\rho}

\def\bbr{{I\!\!R}}
\def\bbn{I\!\!N}
\def\a{\alpha}
\def\b{\beta}
\def\d{\delta}

\def\vr{\vert}

\def\L{\Lambda}

\def\rank{{\rm rank}}

\def\riga{-\kern-4pt - \kern-4pt -}
\font\fat=cmsy10 scaled\magstep5

\def\Bbullet{\raise-3pt\hbox{\fat\char"0F}}

\def\ca{{\cal A}}  \def\cc{{\cal C}}
\def\cd{{\cal D}} \def\ce{{\cal E}} \def\cf{{\cal F}}
\def\cg{{\cal G}} \def\ch{{\cal H}} 
 \def\ck{{\cal K}} 
\def\cm{{\cal M}} \def\cn{{\cal N}} 
\def\cp{{\cal P}}  
 \def\ct{{\cal T}}

\def\Th{\Theta}

\def\ido{intertwining differential operator}
\def\idos{intertwining differential operators}

 \def\ha{{\textstyle{\frac{1}{2}}}}

\def\eps{\epsilon}

\def\ca{{\cal A}}
\def\nn{\nonumber}

\def\nt{\noindent}

\def\lra{\longleftrightarrow}

\input epsf.tex
\newcount\figno
\def\fig#1#2#3{
\par\begingroup\parindent=0pt\leftskip=1cm\rightskip=1cm\parindent=0pt
\baselineskip=11pt \global\advance\figno by 1 
\epsfxsize=#3 \centerline{\epsfbox{#2}} \vskip 12pt
#1\par
\endgroup\par}
\def\figlabel#1{\xdef#1{\the\figno}}
\def\encadremath#1{\vbox{\hrule\hbox{\vrule\kern8pt\vbox{\kern8pt
\hbox{$\displaystyle #1$}\kern8pt} \kern8pt\vrule}\hrule}}

\begin{document}

 \title{Invariant Differential Operators for Non-Compact Lie
Groups:\
the\  SO$^*$(12)\ Case}

\author{V.K. Dobrev}

\address{Institute for Nuclear
Research and Nuclear Energy,
 Bulgarian Academy of Sciences,  72
Tsarigradsko Chaussee, 1784 Sofia, Bulgaria}

\ead{dobrev@inrne.bas.bg}

 \begin{abstract}
In the present paper we continue the project of systematic
construction of invariant differential operators on the example of
the non-compact  algebra  $so^*(12)$.  We give
 the main  multiplets of indecomposable elementary
representations.   Due to the recently established
parabolic relations the multiplet classification results are valid also for the algebra
$so(6,6)$ with suitably chosen maximal parabolic subalgebra.
\end{abstract}

\section{Introduction}

Invariant differential operators   play very important role in the
description of physical symmetries.
In a recent paper \cite{Dobinv} we started the systematic explicit
construction of invariant differential operators. We gave an
explicit description of the building blocks, namely, the parabolic
subgroups and subalgebras from which the necessary representations
are induced. Thus we have set the stage for study of different
non-compact groups.

 In the present paper we  focus on the  algebra  $so^*(12)$.
 The algebras $so^*(4n)$  are form a subclass of the class of algebras,
 which we call 'conformal Lie algebras' in \cite{Dobeseven},
which have very similar properties to the canonical conformal
algebras of  Minkowski space-time.
In our subclass we have the algebras: $so^*(4)$, $so^*(8)$, $so^*(12)$, ...
However the first two cases are reduced to well known conformal algebras due to the coincidences:
~$so^*(4) \cong so(3) \oplus so(2,1)$, ~$so^*(8) \cong so(6,2)$.
Thus, the algebra ~$so^*(12)$~ is the lowest nontrivial member of our subclass.

 This paper is a   sequel of \cite{Dobinv} and \cite{Dobparab}
  and due to the lack of space we refer to these papers
 for motivations and extensive list of
literature on the subject.

\section{Preliminaries}

 Let $G$ be a semisimple non-compact Lie group, and $K$ a
maximal compact subgroup of $G$. Then we have an Iwasawa
decomposition ~$G=KA_0N_0$, where ~$A_0$~ is abelian simply
connected vector subgroup of ~$G$, ~$N_0$~ is a nilpotent simply
connected subgroup of ~$G$~ preserved by the action of ~$A_0$.
Further, let $M_0$ be the centralizer of $A_0$ in $K$. Then the
subgroup ~$P_0 ~=~ M_0 A_0 N_0$~ is a minimal parabolic subgroup of
$G$. A parabolic subgroup ~$P ~=~ M A N$~ is any subgroup of $G$
  which contains a minimal parabolic subgroup.

The importance of the parabolic subgroups comes from the fact that
the representations induced from them generate all (admissible)
irreducible representations of $G$ \cite{Lan,Zhea,KnZu}.

Let ~$\nu$~ be a (non-unitary) character of ~$A$, ~$\nu\in\ca^*$,
let ~$\mu$~ fix an irreducible representation ~$D^\mu$~ of ~$M$~ on
a vector space ~$V_\mu\,$.

 We call the induced
representation ~$\chi =$ Ind$^G_{P}(\mu\otimes\nu \otimes 1)$~ an
~{\it elementary representation} of $G$ \cite{DMPPT}.   Their spaces of functions are:
\eqn{fun} \cc_\chi ~=~ \{ \cf \in C^\infty(G,V_\mu) ~ \vr ~ \cf
(gman) ~=~ e^{-\nu(H)} \cdot D^\mu(m^{-1})\, \cf (g) \} \eeq where
~$a= \exp(H)\in A$, ~$H\in\ca\,$, ~$m\in M$, ~$n\in N$. The
representation action is the $left$ regular action: \eqn{lrr}
(\ct^\chi(g)\cf) (g') ~=~ \cf (g^{-1}g') ~, \quad g,g'\in G\ .\eeq

For our purposes we need to restrict to ~{\it maximal}~ parabolic
subgroups ~$P$, so that $\rank\,A=1$. Thus, for our  representations
the character ~$\nu$~ is parameterized by a real number ~$d$,
called the conformal weight or energy.

An important ingredient in our considerations are the ~{\it
highest/lowest weight representations}~ of ~$\cg$. These can be
realized as (factor-modules of) Verma modules ~$V^\L$~ over
~$\cg^\bbc$, where ~$\L\in (\ch^\bbc)^*$, ~$\ch^\bbc$ is a Cartan
subalgebra of ~$\cg^\bbc$, weight ~$\L = \L(\chi)$~ is determined
uniquely from $\chi$ \cite{Har,Dob}.

Actually, since our ERs will be induced from finite-dimensional
representations of ~$\cm$~ (or their limits)  the Verma modules are
always reducible. Thus, it is more convenient to use ~{\it
generalized Verma modules} ~$\tV^\L$~ such that the role of the
highest/lowest weight vector $v_0$ is taken by the
 space ~$V_\mu\,v_0\,$. For the generalized
Verma modules (GVMs) the reducibility is controlled only by the
value of the conformal weight $d$. Relatedly, for the \idos{} only
the reducibility w.r.t. non-compact roots is essential.

One main ingredient of our approach is as follows. We group the
(reducible) ERs with the same Casimirs in sets called ~{\it
multiplets} \cite{Dobmul,Dob}. The multiplet corresponding to fixed
values of the Casimirs may be depicted as a connected graph, the
vertices of which correspond to the reducible ERs and the lines
between the vertices correspond to intertwining operators. The
explicit parametrization of the multiplets and of their ERs is
important for understanding of the situation.

In fact, the multiplets contain explicitly all the data necessary to
construct the \idos{}. Actually, the data for each \ido{} consists
of the pair ~$(\b,m)$, where $\b$ is a (non-compact) positive root
of ~$\cg^\bbc$, ~$m\in\bbn$, such that the BGG \cite{BGG} Verma
module reducibility condition (for highest weight modules) is
fulfilled: \eqn{bggr} (\L+\r, \b^\vee ) ~=~ m \ , \quad \b^\vee
\equiv 2 \b /(\b,\b) \ .\eeq When \eqref{bggr} holds then the Verma
module with shifted weight ~$V^{\L-m\b}$ (or ~$\tV^{\L-m\b}$ ~ for
GVM and $\b$ non-compact) is embedded in the Verma module ~$V^{\L}$
(or ~$\tV^{\L}$). This embedding is realized by a singular vector
~$v_s$~ determined by a polynomial ~$\cp_{m,\b}(\cg^-)$~ in the
universal enveloping algebra ~$(U(\cg_-))\ v_0\,$, ~$\cg^-$~ is the
subalgebra of ~$\cg^\bbc$ generated by the negative root generators
\cite{Dix}.
More explicitly, \cite{Dob}, ~$v^s_{m,\b} = \cp^m_{\b}\, v_0$ (or ~$v^s_{m,\b} = \cp^m_{\b}\, V_\mu\,v_0$ for GVMs).
  Then there exists \cite{Dob} an \ido{} \eqn{lido}
\cd^m_{\b} ~:~ \cc_{\chi(\L)} ~\llr ~ \cc_{\chi(\L-m\b)} \eeq given
explicitly by: \eqn{mido}\cd^m_{\b} ~=~ \cp^m_{\b}(\widehat{\cg^-})
\eeq where ~$\widehat{\cg^-}$~ denotes the $right$ action on the
functions ~$\cf$, cf. \eqref{fun}.

\section{The non-compact Lie algebras $so^*(12)$}

\nt   The group ~$G=SO^*(2n)$~ consists of all matrices in
$SO(2n,\bbc)$ which commute with a real skew-symmetric matrix times
the complex conjugation operator $C$~: \eqn{SOs} SO^*(2n) \doteq \{\
g\in SO(2n,\bbc) ~|~ J_n C g= g J_n C \} \eeq The Lie algebra
$\cg=so^*(2n)$ is given by: \eqnn{sos} so^*(2n) &\doteq& \{\
X\in so(2n,\bbc) ~|~ J_n C X= X J_n C  \}   = \\ &=&   \{\ X=
\begin{pmatrix}a & b \cr - {\bar b} & {\bar a}\end{pmatrix}
 ~|~\ a,b\in gl(n,\bbc), ~~
^ta = -a, ~~ b^\dag = b \  \}\ .\nn\ee $\dim_R\,\cg = n(2n-1)$,
$\rank\,\cg =n$.

The Cartan involution is given   by: ~$\Th X = -X^\dag$. Thus, $\ck
\cong u(n)$: \eqn{susc} \ck = \{\ X= \begin{pmatrix} a & b \cr -b & a\end{pmatrix} ~|~\
 a,b\in gl(n,\bbc), ~~^ta = -a= -{\bar a}, ~~b^\dag = b = {\bar b}
\ \}\ . \eeq Thus,  $\cg = so^*(2n)$ has discrete series representations
and highest/lowest weight representations.  The complimentary
space $\cp$ is given by: \eqn{susp}\cp = \{\ X= \begin{pmatrix}a & b \cr b
& -a\end{pmatrix} ~|~\ a,b\in gl(n,\bbc) \ , ~~ ^ta = -a = {\bar a}, ~~b^\dag =
b = -{\bar b}\
 \ \}\ . \eeq
$\dim_R\,\cp = n(n-1)$. The split rank is  ~$r\equiv [n/2]$.
The subalgebras ~$\cn_0^\pm$~ which form the root spaces of the root
system $(\cg,\ca_0)$ are of real dimension ~$n(n-1) - [n/2]$.

The maximal parabolic subalgebras have $\cm$-factors as follows \cite{Dobinv}:
\eqn{maxsose} \cm^{\rm max}_j ~=~ so^*(2n-4j) \oplus su^*(2j)  \ , ~~~j=1,\ldots,r\ . \eeq
The ~$\cn^\pm$~  factors in the maximal parabolic subalgebras
have dimensions:  ~$\dim\,(\cn^\pm_j)^{\rm max} ~=$ $j(4n-6j-1)$.

For even ~$n=2r$~ we choose a ~{\it maximal} parabolic ~$\cp=\cm\ca\cn$~ such that
~$\ca\cong so(1,1)$, ~$\cm ~=~ \cm^{\rm max}_r ~=~ ~su^*(n)$. We note also that
\eqn{colial} \ck^\bbc \cong
u(1)^\bbc \oplus sl(n,\bbc)  \cong \ca^\bbc \oplus \cm^\bbc \eeq
Thus,  the factor ~$\cm$~ has the same
finite-dimensional (nonunitary) representations as the
finite-dimensional (unitary) representations of  the semi-simple
subalgebra of   ~$\ck$. The property \eqref{colial} distinguishes
the class we called 'conformal Lie algebras' \cite{Dobeseven}, to which class
the algebras ~$so^*(4r)$~ belong.

\bigskip

 Further we restrict to our case of study ~$\cg ~=~ so^*(12)$.

\bigskip

We label   the signature of the ERs of $\cg$   as follows:
\eqn{sgnd}  \chi ~=~ \{\, n_1\,, n_{2}\,, n_{3}\,, n_{4}\,, n_{5}\,;\, c\, \} \ ,
\qquad n_j \in \bbz_+\ , \quad c =
d- \ftha \eeq where the last entry of ~$\chi$~ labels the characters of
$\ca\,$, and the first five entries are labels of the
finite-dimensional (nonunitary) irreps of $\cm=su^*(6)$ when all $n_j>0$  or
limits of the latter when some $n_j=0$.

Below we shall use the following conjugation on the
finite-dimensional entries of the signature: \eqn{conu}
(n_1,n_{2},n_3,n_{4},n_{5})^* ~\doteq~
(n_{5},n_{4},n_3,n_2,n_{1})  \ . \eeq

The ERs in the multiplet are related also by intertwining integral
  operators  introduced in \cite{KnSt}. These operators are defined
for any ER,   the general action
being: \eqnn{knast}  && G_{KS} ~:~ \cc_\chi ~ \llr ~ \cc_{\chi'} \
,\cr &&\chi ~=~ \{\, n_1,\ldots ,n_{5}
\,;\, c\, \} \ , \qquad \chi' ~=~ \{\,
(n_1,\ldots,n_{5})^* \,;\, -c\, \} . \ee

Further, we need the root system of the complexification
~$\cg^\bbc = so(12,\bbc)$~.  The positive roots are given standardly as:
\eqnn{sunnpos} \a_{ij} ~&=&~ \eps_i -  \eps_{j}\ , \quad 1 \leq i < j \leq
6 \ , \cr \b_{ij} ~&=&~ \eps_i +  \eps_{j}  \ , \quad 1 \leq i<j \leq 6 \ee
where ~$\eps_i$~  are standard orthonormal basis: ~$\lg \eps_i,\eps_j\rg = \d_{ij}\,$.
The compact roots are ~$\a_{ij}$ - they  form (by
restriction) the root system of the semisimple part of ~$\ck^\bbc \cong
su^*(6)^\bbc$,  while the roots ~$\b_{ij}$~ are noncompact.

Further, we give the correspondence between the signatures $\chi$
and the highest weight $\L$. The connection is through the Dynkin
labels:    \eqn{dynk} m_i ~\equiv~ (\L+\r,\a^\vee_i) ~=~ (\L+\r,
\a_i )\ , \quad i=1,\ldots,6, \eeq where ~$\L = \L(\chi)$, ~$\r$ is
half the sum of the positive roots of ~$\cg^\bbc$.  The explicit
connection is: \eqn{rela} n_i = m_i \ , \quad  c
 ~=~  -\,\ha(   m_1+ 2m_{2} + 3m_3 + 4m_{4} + 2m_{5}+ 3 m_6)
 \eeq

Finally, we remind that according to  \cite{Dobparab} the above considerations
  are applicable also for the algebra ~$so(6,6)$~
with parabolic ~$\cm$-factor ~$sl(6,\bbr)$.

\section{Main multiplets of SO$^*$(12)}

The main multiplets   are in 1-to-1 correspondence
with the finite-dimensional irreps of ~$so^*(12)$, i.e., they are
labelled by  the  six  positive Dynkin labels    ~$m_i\in\bbn$.

The number of ERs/GVMs in the corresponding multiplets is  (cf. \cite{Dobparab}):
\eqn{weyln}
\vr W(\cg^\bbc,\ch^\bbc)\vr\, /\, \vr W(\ck^\bbc,\ch^\bbc)\vr \ =\
\vr W(so(12,\bbc))\vr\, /\, \vr W(sl(6,\bbc))\vr\ =\ 32  \eeq
where ~$\ch$~ is a Cartan subalgebra of both ~$\cg$~ and ~$\ck$.

The signatures of the  32 ERs/GVMs in a  main multiplet are
given in the following pair-wise manner:
\eqnn{tabltri}  \chi_0^\pm
   &=&    \{  (
 m_1,
 m_2,
 m_3
 m_4,
 m_5)^\pm ;\  \pm  \ha(m_1+2m_2+3m_3+4m_4+2m_5+3m_6)  \} \nn\\
\chi_a^\pm      &=&      \{  (
 m_1,
 m_2,
 m_{3},
 m_{4,6},
 m_5)^\pm ;\  \pm \ha (m_1+2m_2+3m_3+4m_4+2m_5+m_6)   \} \nn\\
\chi_b^\pm    &=&    \{  (
 m_1, m_{2},
 m_{34},
 m_{6},
 m_{45})^\pm ;\  \pm \ha (m_1+2m_2+3m_3+2m_4+2m_5+m_6)   \} \nn\\
\chi_c^\pm      &=&      \{  (
m_{1}, m_{2},
 m_{35},
 m_{6},
 m_4)^\pm ;\  \pm \ha (m_1+2m_2+3m_3+2m_4 +m_6)   \} \nn\\
\chi_{c'}^\pm          &=&          \{  (
 m_{1},
 m_{23},
 m_{4},
 m_{6},
 m_{35})^\pm ;\  \pm \ha (m_1+2m_2+m_3+2m_4+2m_5+m_6)    \} \nn\\
\chi_d^\pm        &=&        \{  (
m_{1},
 m_{23},
 m_{45},
 m_{6},
 m_{34})^\pm ;\  \pm \ha (m_1+2m_2+m_3+2m_4 +m_6)     \} \nn\\
\chi_{d'}^\pm        &=&        \{  (
 m_{12},
 m_{3},
 m_{4},
 m_{6},
 m_{25})^\pm ;\  \pm  \ha (m_1 +m_3+2m_4+2m_5+m_6)  \} \nn\\
\chi_e^\pm      &=&      \{  (
 m_{1},
 m_{24},
 m_{5},
 m_{4,6},
 m_{3})^\pm ;\  \pm  \ha (m_1+2m_2+m_3  +m_6)     \}\nn\\
\chi_{e'}^\pm        &=&        \{  (
 m_{12},
 m_{3},
 m_{45},
 m_{6},
 m_{24})^\pm ;\  \pm  \ha (m_1 +m_3+2m_4 +m_6)  \} \nn\\
 \chi_{e''}^\pm        &=&        \{  (
 m_{2},
 m_{3},
 m_{4},
 m_{6},
 m_{15})^\pm ;\  \pm  \ha (-m_1 +m_3+2m_4+2m_5+m_6)  \} \nn\\
 \chi_{f}^\pm        &=&        \{  (
 m_{12},
 m_{34},
 m_{5},
 m_{4,6},
 m_{23})^\pm ;\  \pm  \ha (m_1 +m_3 +m_6)  \} \nn\\
 \chi_{f'}^\pm        &=&        \{  (
 m_{2},
 m_{3},
 m_{45},
 m_{6},
 m_{14})^\pm ;\  \pm  \ha (-m_1 +m_3+2m_4 +m_6)  \} \nn\\
\chi_{f''}^\pm      &=&      \{  (
 m_{1},
 m_{24,6},
  m_{5}
 m_{4},
 m_{3})^\pm ;\  \pm  \ha (m_1+2m_2+m_3 - m_6)  \} \nn\\
 \chi_{g}^\pm        &=&        \{  (
 m_{13},
 m_{4},
 m_{5},
 m_{34,6},
 m_{2})^\pm ;\  \pm  \ha (m_1 -m_3 +m_6)  \} \nn\\
 \chi_{g'}^\pm        &=&        \{  (
 m_{12},
 m_{34,6},
 m_{5},
 m_{4},
 m_{23})^\pm ;\  \pm  \ha (m_1 +m_3 -m_6)  \} \nn\\
 \chi_{g''}^\pm        &=&        \{  (
 m_{2},
 m_{34},
 m_{5},
 m_{4,6},
 m_{13})^\pm ;\  \pm  \ha (-m_1 +m_3 +m_6)  \}
     \ee
 where ~$(k_1,k_2,k_3,k_4,k_5)^- = (k_1,k_2,k_3,k_4,k_5)$,
~$(k_1,k_2,k_3,k_4,k_5)^+ = (k_1,k_2,k_3,k_4,k_5)^*$.
 They are given explicitly in Fig.~1.
The pairs ~$\L^\pm$~ are symmetric w.r.t. to the bullet in the
middle of the figure - this represents the Weyl symmetry realized by
the Knapp-Stein operators \eqref{knast}:~ $G_{KS} ~:~ \cc_{\chi^\mp}
\lra \cc_{\chi^\pm}\,$.

Matters are arranged so that in every multiplet only the ER with
signature ~$\chi_0^-$~ contains a finite-dimensional nonunitary
subrepresentation in  a finite-dimensional subspace ~$\ce$. The
latter corresponds to the finite-dimensional   irrep of ~$so^*(12)$~ with
signature ~$\{ m_1\,, \ldots\,, m_6 \}$.   The subspace ~$\ce$~ is annihilated by the
operator ~$G^+\,$,\ and is the image of the operator ~$G^-\,$. The
subspace ~$\ce$~ is annihilated also by the \ido{} acting from
~$\chi^-_0$~ to ~$\chi^-_a\,$.
 When all ~$m_i=1$~ then ~$\dim\,\ce = 1$, and in that case
~$\ce$~ is also the trivial one-dimensional UIR of the whole algebra
~$\cg$. Furthermore in that case the conformal weight is zero:
~$d=\ftha+c=\ftha-\ha(m_1+2m_2+3m_3+4m_4+2m_5+3m_6)_{\vert_{m_i=1}}=0$.

In the conjugate ER ~$\chi_0^+$~ there is a unitary discrete series
subrepresentation in  an infinite-dimen\-sional subspace  $\cd$. It
is annihilated by the operator  $G^- $,\ and is the image of the
operator  $G^+ $.

 All the above is valid also for the algebra ~$so(6,6)$, cf.  \cite{Dobparab}.
 However, the latter algebra does not have highest/lowest weight
   representations while the algebra
 ~$so^*(12)$~ has  highest/lowest weight series representations.

Thus, for  ~$so^*(12)$~  the   ER with signature
~$\chi_0^+$~ contains both a holomorphic discrete series representation
and a conjugate
anti-holomorphic discrete series representation. The direct sum of the holomorphic
and the antiholomorphic representations spaces form the
invariant subspace ~$\cd$~ mentioned above.
Note that the corresponding lowest weight GVM
is infinitesimally equivalent only to the holomorphic discrete
series, while the conjugate highest weight GVM is infinitesimally
equivalent to the anti-holomorphic discrete series.

In Fig.~1   we use the notation: ~$\L^\pm = \L(\chi^\pm)$.
Each \ido\ is represented by an
arrow accompanied by a symbol ~$i_{jk}$~ encoding the root
~$\a_{jk}$~ and the number $m_{\a_{jk}}$ which is involved in
the BGG criterion.  This notation is used to save space, but it can
be used due to the fact that only \idos\ which are
non-composite are displayed, and that the data ~$\b,m_\b\,$, which
is involved in the embedding ~$V^\L \lra V^{\L-m_\b,\b}$~ turns out
to involve only the ~$m_i$~ corresponding to simple roots, i.e., for
each $\b,m_\b$ there exists ~$i = i(\b,m_\b,\L)\in \{
1,\ldots,6\}$, such that ~$m_\b=m_i\,$. Hence in Fig. 1. the data
~$\a_{jk}\,$,~$m_{\a_{jk}}$~ is represented by ~$i_{jk}$~
on the arrows.

\fig{}{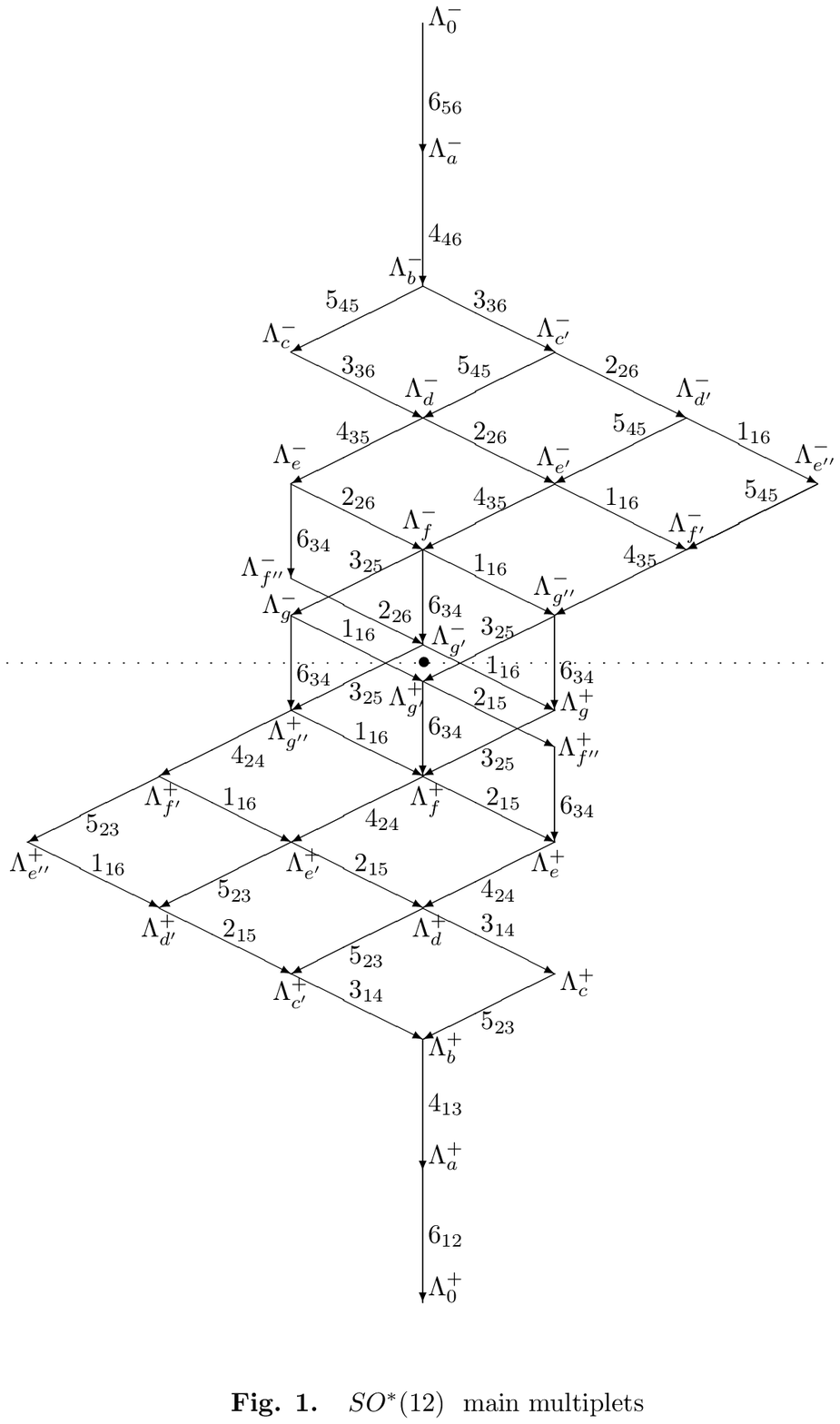}{12cm}

\np

\section*{Acknowledgments.}\nt
 The author would like to thank the Organizers for the kind
hospitality and invitation to present a talk at the XXX International Colloquium
on Group Theoretical Methods in Physics (Ghent, July 2014).

\end{document}